\documentclass{elsart}

\usepackage{xy}
\xyoption{all}

\usepackage{amsmath}
\usepackage{amssymb}
\usepackage{theorem}
\newfont{\cyr}{wncyr10 at 11 pt}

\newcommand {\art}[6]{{\sc #1:} {#2.} {\em #3} {\bf #4} {(#5),} {#6.}}
\newcommand {\book}[5]{{\sc #1:} {``#2."} {#3,} {#4} {(#5).}}
\newcommand {\samp}[8]{{\sc #1:} {#2,} {{\em in:} ``#3,"} {(#4),}
                                  {#5,} {#6,} {(#7),} {#8.}}

\newcommand {\toappear}[3]{{\sc #1:} {#2.} {\em #3,} {\sl{to appear}.}}
\newcommand {\submitted}[2]{{\sc #1:} {#2.} {\sl submitted}.}

\newcommand {\preprint}[2]{{\sc #1:} {#2,} {\sl preprint.}}
\newcommand {\misc}[3]{{\sc #1:} {#2,} {\sl #3.}}

\renewcommand \phi {\varphi}

\newcommand \Yu {\mbox{\cyr \symbol{16}}}
\newcommand {\YU}[1]{\Yu{(}#1{)}}
\newcommand {\YUp}[1]{{\Yu}^{+}{(}#1{)}}
\newcommand {\YUm}[1]{{\Yu}^{-}{(}#1{)}}

\newcommand \parto {(\omega)^\omega}
\newcommand \parti {(\omega)}
\newcommand \parlo {(\omega)^{<\omega}}
\newcommand \parce {(\omega)_{{}^{\ceq}}}
\newcommand \parfe {(\omega)_{{}^{\feq}}}
\newcommand \partoce {(\omega)^{\omega}_{{}^{\ceq}}}

\newcommand \reals {[\omega ]^\omega }

\newcommand \mo {{\it{m.o.}}}

\newcommand \ort {\bot}

\newcommand \und {\wedge}

\newcommand \mmin {\operatorname{min}}

\newcommand \mMin {\operatorname{Min}}

\newcommand \feq {\sqsupseteq}
\newcommand \ceq {\sqsubseteq}

\newcommand \kap {\sqcap}

\newcommand \kup {\sqcup}
\newcommand \seg {\preccurlyeq}

\newcommand \cF {{\mathcal F}}
\newcommand \cP {{\mathcal P}}

\newcommand \cB {{\mathcal B}}
\newcommand \cO {{\mathcal O}}
\newcommand \cA {{\mathcal A}}

\newcommand \fS {{\mathfrak S}}
\newcommand \fH {{\mathfrak H}}

\newcommand \fh {{\mathfrak h}}

\newcommand \st {{\scriptscriptstyle{\bigstar}}}

\newcommand \bfN {{\mathbf N}}

\newcommand \bfV {{\mathbf V}}

\newcommand \bbL {{\mathbb L}}
\newcommand \bbU {{\mathbb U}}

\newcommand \CH {{\sf CH}}

\newcommand \ZFC {{\sf ZFC}}

\newcommand \CS {{\mathbf{Set}}}
\newcommand{\Sub}[1]{{{\rm Sub}_\CS(#1)}}
\newcommand{\Cosub}[1]{{{\rm Cosub}_\CS(#1)}}

\newcommand \la {\langle}
\newcommand \ra {\rangle}
\newcommand \subs {\subseteq}

{\theorembodyfont{\rmfamily}
\newtheorem {nummer}{ }[section]

\newtheorem {thrm}[nummer]{Theorem}

\newtheorem {propn}[nummer]{Proposition}

\newtheorem {lm}[nummer]{Lemma}
\newtheorem {fct}[nummer]{Fact}
\newtheorem {obs}[nummer]{Observation}
}

\def\qed{{\unskip\nobreak\hfil\penalty50\hskip8mm\hbox{}
  \nobreak\hfil
  {{\small{\bf{q.e.d.}}}}\parfillskip=0mm \par\smallskip}}

\begin{document}
\begin{frontmatter}
\title{Ultrafilter Spaces on the Semilattice of Partitions}
\author{Lorenz Halbeisen}
\address{
Department of
Mathematics, University of California, Evans Hall 938, Berkeley, CA 94720,
USA, {\tt halbeis@math.berkeley.edu}}
\author{Benedikt L\"owe}
\address{
Mathematisches Institut, Rheinische
Friedrich--Wilhelms--Universit\"at Bonn, Beringstra{\ss}e 6, 53115
Bonn, Germany, {\tt loewe@math.uni-bonn.de}}
\thanks{The first author wishes to thank the
Swiss National Science Foundation for supporting him, and the second
author wishes to thank the Studienstiftung des deutschen Volkes and
the Graduiertenkolleg ``Algebraische, analytische und geo\-metrische
Methoden und ihre Wechselwirkung in der modernen Mathematik'' in Bonn
for travel grants. Both authors would like to extend their gratitude
to Andreas Blass (Ann Arbor MI) for sharing his views on dualization
with them. }
\begin{abstract}
The Stone-\v{C}ech compactification of the natural numbers
$\beta\omega$ (or equi\-valently, the space of ultrafilters on the
subsets of $\omega$) is a well-studied space with interesting
properties. Replacing the subsets of $\omega$
by partitions of
$\omega$ in the construction of the
ultrafilter space gives non-homeomorphic spaces of
partition ultrafilters corresponding to $\beta\omega$.
We develop a general framework for spaces of this type and
show that the spaces of partition ultrafilters still have
some of the nice properties of $\beta\omega$, even though none of them is
homeomorphic to $\beta\omega$. Further, in a particular space, the
minimal height of a tree $\pi$-base and $P$-points are investigated.\\[0,3cm]
{\it Keywords:} Compactness; Cosubobjects; Partitions; P-Points;
Separability; Shattering families; Stone topologies on semilattices;
Ultrafilter spaces\\[0,3cm]
2000 {\it Mathematics Subject Classification}: {\bf 54A05 54D80} 54A10 03E17 03E35
\end{abstract}
\end{frontmatter}

\section{Introduction}\label{sec:int}

This paper is a glance at a generalization of the space of
ultrafilters over $\omega$, the Stone-\v{C}ech compactification of
the natural numbers, or just $\beta\omega$. The space $\beta\omega$
and its remainder, the space $\beta\omega\setminus\omega$, are
well-studied spaces with a lot of interesting properties. For example
both spaces are quasicompact and Hausdorff and therefore,
$C(\beta\omega)$ and $C(\beta\omega\setminus\omega)$ are Banach
algebras.

\bigskip

We shall provide the reader with a very general approach to
ultrafilter spaces on arbitrary semilattices, but the main focus of
this paper is a particular class of ultrafilter spaces because of
their intimate category theoretical connection with the
Stone-\v{C}ech compactification: spaces of ultrafilters on the
semilattice of partitions.

These objects are indeed the dualization of $\beta\omega$ in the
following category theoretical sense:

\smallskip

Look at the category of sets $\CS$. As usual for any object $M\in{\rm
Obj}_\CS$ we define the subobjects of $M$ to be the equivalence
classes of monos\footnote{Two monos $f:A\to M$ and $g:B\to M$ are
defined to be equivalent iff there is an isomorphism $h:A\to B$ such
that $f=gh$.} ({\sl i.e.}, in the category of sets, injections)
with codomain $M$. The
collection $\Sub{M}$ of subobjects of $M$ is partially ordered by
$$[f]\leq [g]
:\iff \exists h:A\to B (f = gh),$$ when $f:A\to M$ and $g:B\to M$.

Indeed, in this case $\Sub{M}$ is a Boolean algebra with greatest
element $[{\rm id}_M]$ and least element $[o_M]$ where $o_M$ is the
unique morphism with codomain $M$ whose domain is the initial object
of the category $\CS$ (the empty set).

The Stone-\v{C}ech compactification $\beta\omega$ is the space of
ultrafilters in the Boolean algebra $\Sub{\omega}$.

\medskip

We can now dualize by reversing all occurring arrows: The dualization
of $\Sub{M}$ is the collection $\Cosub{M}$ of all equivalence classes
of epis\footnote{Two epis $f:M\to A$ and $g:M\to B$ are defined to be
equivalent iff there is an isomorphism $h:B\to A$ such that $f=hg$.}
({\sl i.e.}, in the category of sets, surjections)
with domain $M$, which we will call the
{\bf cosubobjects} of $M$. Again, this collection is partially
ordered by $$[f]\leq [g] :\iff \exists h:B\to A (f = hg),$$ when
$f:M\to A$ and $g:M\to B$.

\medskip

Moving from the abstract to the concrete, in the category of sets, a
cosubobject of $\omega$ is just an equivalence class of surjective
functions with domain $\omega$ modulo permutation of their ranges.
This gives us a partition of $\omega$ by looking at the preimages of
singletons of elements of the range of such functions. Note that
following this translation, the relation $\leq$ defined on
$\Cosub{M}$ just gives us the ``is coarser than'' relation on
partitions.

Thus $\la \Cosub{\omega},\leq\ra$ is again a partially ordered
structure and its space of ultrafilters is in this sense the
dualization of the Stone-\v{C}ech compactification $\beta\omega$.

\smallskip

The important distinction between infinite sets and finite sets in
the Stone-\v{C}ech compactification that allows to distinguish
between principal ultrafilters ({\sl i.e.}, the representants of the
countable dense subset $\omega$ of $\beta\omega$) and non-principal
ultrafilters becomes dualized to the distinction between partitions
into infinitely many classes and partitions into finitely many
classes.\footnote{Note that the notion of finiteness of a subobject
of $\omega$ can be categorially expressed by use of the Dedekind
formalization of finiteness (using the Axiom of Choice) : $[f]$ is
finite if and only if every mono from the domain of $f$ to itself is
an isomorphism.}

\smallskip

But the move from subobjects to cosubobjects changes quite a lot: We
will see in this paper that $\la \Cosub{\omega},\leq\ra$ cannot be a
Boolean algebra. An immediate consequence of this lack of a
complementation function is that some variations of techniques, that
were merely different viewpoints in the case of the Stone-{\v C}ech
compactification, now actually give different topological spaces. We
will show that spaces that are homeomorphic in the classical case
fall apart in the dual case, especially two consequences of
compactness, countable compactness and the Hausdorff separation
property, belong to two different spaces in the dual case --- so none
of the possible dualization of the Stone-\v{C}ech compactification is
compact anymore.

This alone should be enough motivation to delve deeper into that
subject matter to get more information about these spaces and find
the most natural dualization of $\beta\omega$. We close our paper
with an extensive list of projects and open problems that result from
these non-homeomorphicity results. Many areas of application for the
Stone-\v{C}ech compactification that are nowadays very well
understood deserve to be explored in our dual case.

\bigskip

We will restrict our attention in this paper to partitions which
consist solely of infinite blocks. There is no innate category
theoretical reason behind this, but we believe that the additional
information that the size of a block might convey could add unwanted
combinatorial phenomena to the theory of spaces of partition
ultrafilters. After all, set-theoretically speaking, the elements of
a subset of $\omega$ also do not carry an additional information
about their size, so we try to avoid hidden information by
restricting our attention to partitions with large blocks, and thus
receive some sort of homogeneity.

For readers interested in other approaches, we mention this
restriction (and a possible lifting of it) in Section~\ref{sec:open}.

\section{Ultrafilter Spaces on Semilattices}\label{sec:UltSem}

In this section we define topologies on the set of ultrafilters on
semilattices.

\subsection{Semilattices and Partitions}

A {\bf semilattice} $\bbL=\la L,\seg,{\mathbf 0}\ra$ consists of a
set $L$, a least element ${\mathbf 0}$ and a partial ordering
$\seg$ on $L$ such that for all $x,y\in L$ we have the following:
There is a $z\in L$ with $z\seg x$ and $z\seg y$, and for every $w$
with $w\seg x$ and $w\seg y$ we have $w\seg z$ (the {\bf infimum} of
$x$ and $y$) which we denote as usual by $x\und y$. Furthermore, for each $x\in L$,
the least
element ${\mathbf 0}\in L$ should satisfy ${\mathbf 0}\und x={\mathbf 0}$.
A semilattice without
a least element can easily be supplemented by one.

Furthermore, an element $x\in L$ is called an {\bf atom} if for all
$y\seg x$ we have either $y=x$ or $y={\mathbf 0}$. A semilattice is
said to be {\bf downward splitting} if below each $x\in L$ which is
not an atom there are $y_0\in L$ and $y_1\in L$ such that ${\mathbf
0}\neq y_0\seg x$, ${\mathbf 0}\neq y_1\seg x$, and $y_0\und y_1 =
{\bf 0}$.

\medskip

Let $\bbL=\la L,\seg\ra$ be a semilattice. Two elements $x,y\in L$
are called {\bf ortho\-gonal}, and we write $x\ort y$ if $x\und y
={\mathbf 0}$. Otherwise, they are called {\bf compatible}. If we
want to stress the connection between the relations $\ort$ and
$\seg$, we write $\ort_\seg$.

A semilattice $\bbL$ is called {\bf complemented} if there is a
function ${\sim}
: L\to L$ satisfying
\begin{itemize}
\item[(C1)] $\forall x\in L ({\sim}{\sim} x = x)$,
\item[(C2)] $\forall x,y\in L(y\und x = {\mathbf 0}\leftrightarrow
y\seg{\sim} x)$.
\end{itemize}

Complemented semilattices are extremely well-behaved: We can define
the reverse relation by stipulating $x \succcurlyeq y$ iff
$({\sim}x\seg {\sim}y)$. Then ${\bbL}^\updownarrow := \la
L,\succcurlyeq\ra$ is a semilattice with least element
${\sim}{\mathbf 0}$ (where ${\mathbf 0}$ is the least element of
$\bbL$) and the semilattice ${\bbL}^\updownarrow$ is isomorphic to
${\bbL}$ via the map ${\sim}$.

\bigskip

The semilattice we are mainly interested in is the semilattice of
partitions of $\omega$.

A {\bf partition} $X$ (of $\omega$) consisting of pairwise disjoint,
non-empty sets such that $\bigcup X = \omega$. The elements of a
partition are called the {\bf blocks}.

We elaborated in Section \ref{sec:int} on the possibility of a
categorial definition as cosubobjects of $\omega$ and why we will
only consider partitions of $\omega$ all of whose blocks are infinite
sets. So, in the following the word ``partition'' by convention
always refers to partitions of
$\omega$ all of whose blocks are infinite. We also consider finite
partitions, this means partitions containing finitely many blocks,
and the partition containing only one block is denoted by
$\{\omega\}$. The set of all partitions is denoted by $\parti$, the
set of all partitions containing infinitely (resp.\,finitely) many
blocks is denoted by $\parto$ (resp.\,$\parlo$).

Let $X$ and $Y$ be two partitions. We say $X$ is {\bf coarser} than
$Y$, or that $Y$ is {\bf finer} than $X$ (and write $X\ceq Y$) if
each block of $X$ is the union of blocks of $Y$. Let $X\kap Y$ denote
the finest partition which is coarser than $X$ and $Y$. Similarly,
$X\kup Y$ denotes the coarsest partition which is finer than $X$ and
$Y$.

In the following we investigate the semilattices $\parce
:=\la\parti,\ceq\ra$ and $\parfe :=\la\parti\cup\{{\bf 0}\},\feq\ra$,
and in Section~\ref{sec:inf} we will investigate $\partoce
:=\la\parto,\ceq\ra$. Notice that the least element ${\mathbf 0}$ in
$\parce$ is $\{\omega\}\in\parti$, whereas the set $\parfe$ does not
have a least element on its own, so we have to add
one.\footnote{Every partition can be properly refined because all
blocks are infinite, so there is no finest partition.}

\bigskip

Because they figure prominently in the lattice theoretical
description of $\beta\omega$, we also mention the two well-known
semilattices $\cP(\omega)_{{}^\subs}:=\la\cP(\omega),\subs\ra$, where
$\cP(\omega)$ is the power-set of $\omega$, and
$[\omega]^{\omega}_{{}^\subs}:=\la\reals,\subs\ra$, where $\reals$ is
the set of all infinite subsets of $\omega$. As we noted in Section
\ref{sec:int}, $\cP(\omega)_{{}^\subs}$ is not just a semilattice but
a Boolean algebra.

\subsection{Ultrafilters on semilattices}

Let $\bbL = \la L, \seg, \mathbf{0}\ra$ be an arbitrary semilattice.

A family $\cB\subs L$ is called a {\bf filter
base on} $\bbL$ if the following holds: For any $x,y\in\cB$ we have
$x\und y\in\cB$, and ${\mathbf 0}\notin\cB$. If $\cB$ is a filter
base, we shall call $[\cB]:= \{ y : \exists x\in\cB(x\seg y)\}$ the
{\bf filter generated by} $\cB$. A filter base $\cF$ is called a {\bf
filter} if $[\cF] = \cF$. A filter $\cF$ is called an {\bf
ultrafilter} if $\cF$ is not properly contained in any other filter
on ${\bbL}$. A filter $\cF$ is {\bf principal} if there is an $x\in
L$ such that $\cF=[\{x\}]$, otherwise it is called {\bf
non-principal}. As easy consequences of the definition of
ultrafilters (and, in the case of Fact \ref{fct:zorn}, Zorn's Lemma)
we get the following facts:

\begin{fct}\label{fct:equivUlt}
$\cF$ is an ultrafilter on ${\bbL}$ if and only if for any $x \in L$
either $x\in\cF$ or there is a $y\in\cF$ such that $y\und x= {\mathbf
0}$.
\end{fct}

\begin{fct}\label{fct:zorn}
If $X$ is a family of elements of $L$ with the {\bf finite
intersection property} ({\sl i.e.}, for any finite subfamily
$\{x_0,...,x_n\}\subseteq X$ we have $x_0\und ...\und x_n\neq{\bf
0}$), then there is an ultrafilter $\cF$ on $L$ with $X\subseteq
\cF$.
\end{fct}

Let $\YU{\bbL}$ denote the {\bf set of all ultrafilters on} ${\bbL}$.
The Cyrillic letter ``$\Yu$'' for the sound ``yu'' should remind the
reader of the ``u'' in ``ultrafilter''. Note that we make use of the
assumption that the semilattice has a least element: Although we
could get rid of the mention of ${\bf 0}$ in the definition of filter
by postulating $\cF\neq L$ instead of ${\bf 0}\notin\cF$, we cannot
prove Fact \ref{fct:zorn} without the least element. To see this,
look at an arbitrary linear order $\bbL = \la L,\seg\ra$ without
least element. Filters are just endsegments of $\bbL$, but there can
be no maximal proper endsegment. Thus, on this semilattice, there is
no ultrafilter at all.

\subsection{Topologies on ${\YU{\bbL}}$}

We can define topologies on $\YU{\bbL}$ in two different ways:

First define for each $x\in L$ two sets $(x)^+ := \{ p \in\YU{\bbL}
: x \in p\}$ and $(x)^- := \{ p \in\YU{\bbL}
: x \notin p\}=\YU{\bbL}\setminus (x)^+$. Set ${\cO}^+ := \{ (x)^+ :
x \in L\}$ and ${\cO}^- := \{ (x)^- : x\in L\}$ and call the topology
generated by ${\cO}^+$ the {\bf positive topology} $\tau^+$ and the
topology generated by ${\cO}^-$ the {\bf negative topology} $\tau^-$
on $\YU{\bbL}$. (Note that ${\cO}^+$ is a base for $\tau^+$, but
${\cO}^-$ is not necessarily a base for $\tau^-$. This difference
accounts for some of the asymmetries.)

In the following we shall use the notation $\YUp{\bbL}:=
\la\YU{\bbL}, \tau^+\ra$ and $\YUm{\bbL} := \la\YU{\bbL}, \tau^-\ra$.

An immediate consequence of Fact~\ref{fct:equivUlt} is that
$\tau^-\subseteq\tau^+$, since $$(x)^- = \bigcup\{(y)^+ : y\und x =
{\mathbf 0}\}.$$

In the case of complemented semilattices, these two topologies coincide:
To see this, just note that
(C2) implies that ultrafilters contain either $x$ or
${\sim}x$ for each $x\in L$, and that (C1) implies that ${\sim}$ is a
surjective function. Thus, if ${\bbL}$ is complemented, then for each
basic open set $O\in{\cO}^+$ there is an open set
$\tilde{O}\in{\cO}^-$ such that $O =\tilde{O}$, whence $\YUp{\bbL} =
\YUm{\bbL}$.

\bigskip

It is easy to see that $\YUp{\cP(\omega)_{{}^\subs}}$ is just
$\beta\omega$ and that $\YUp{[\omega]^{\omega}_{{}^\subs}}$
is homeomorphic to $\beta\omega\setminus\omega$. Further, since both
semilattices are complemented (for $[\omega]^{\omega}_{{}^\subs}$,
just take the complement if it's infinite and ${\mathbf 0}$ if the
set is cofinite), we get that $\YUp{\cP(\omega)_{{}^\subs}}$ is
homeomorphic to each of the spaces $\YUm{\cP(\omega)_{{}^\subs}}$,
$\YUp{\cP(\omega)_{{}^\supseteq}}$ and
$\YUm{\cP(\omega)_{{}^\supseteq}}$, and that
$\YUp{[\omega]^{\omega}_{{}^\subs}}$ is homeomorphic to
$\YUm{[\omega]^{\omega}_{{}^\subs}}$,
$\YUp{[\omega]^{\omega}_{{}^\supseteq}}$ and
$\YUm{[\omega]^{\omega}_{{}^\supseteq}}$.

We shall call a topological space {\bf principal} if it contains an
open set with just one element. (Proposition \ref{prop:principal}
will explain the choice of the name ``principal'' for this property.)
Being principal is obviously a property preserved under
homeomorphisms, so it is a topological invariant.

\begin{propn}\label{prop:principal}
Let $\bbL$ be a semilattice which splits downward.
Then the following are equivalent:
\begin{enumerate}
\item $\YUp{\bbL}$ is a principal space, and
\item $\YU{\bbL}$ contains a principal ultrafilter.
\end{enumerate}
\end{propn}

\begin{pf}
``(ii) $\Rightarrow$ (i)'': Let $p:=[\{x\}]\in\YU{\bbL}$. Take any
$q\in\YU{\bbL}$ with $x\in q$. Then $q\supseteq p$ and hence by
maximality of $q$ and $p$ (both are ultrafilters), we have $p=q$.
Thus we have $(x)^+ = \{p\}$ and this is our open set with one
element.

``(i) $\Rightarrow$ (ii)'': Now let $\{p\}$ be an open set with one element.
Obviously, such a set must be a basic open set. Let $(x)^+ = \{p\}$.

{\bf Case I:} If $x$ is an atom then $q:=[\{x\}]$ is an ultrafilter
and with $x\in p$ we have $q\subseteq p$. Since both $p$ and $q$ are
ultrafilters, we have $p=q$. Thus $p$ is principal and we are done.

{\bf Case II:} If $x$ is not an atom, we can (by the property of
downward splitting) pick elements ${\mathbf 0}\neq y_0\seg x$ and
${\mathbf 0}\neq y_1\seg x$ such that $y_0\und y_1 = {\mathbf 0}$.
By Fact \ref{fct:zorn},
the sets $\{x,y_0\}$ and $\{x,y_1\}$ can be extended to ultrafilters
$p_1$ and $p_2$. Obviously, both $p_1$ and $p_2$ are elements of
$(x)^+$ and $p_1\neq p_2$, contradicting the assumption that $(x)^+$
is a singleton.\qed\end{pf}

First of all, note that you can't drop the assumption of downward
splitting: Take any dense linear order $\bbL=\la
L,\seg\ra$ with a least element ${\mathbf 0}$.
Then $\YU{\bbL}$ contains just one element (the
ultrafilter of all elements $x\neq{\bf 0}$) and this element is
non-principal, but the space $\YUp{\bbL}$ is principal (since it is a
point).

The nice characterization of Proposition \ref{prop:principal} does
not work in the case of the negative topologies, since the existence
of closed singletons (which would be the analogue of being a principal space
for the negative
topologies) is provable in general regardless of the existence of
principal ultrafilters:

\begin{fct} For any semilattice $\bbL$,
the spaces $\YUp{\bbL}$ and $\YUm{\bbL}$ are ${\rm T}_1$ spaces ({\sl
i.e.}, all singletons are closed).
\end{fct}

\begin{pf}
For any singleton $\{p\}$ look at $\bigcup_{x\notin p} (x)^+$ for the
positive topology and $\bigcup_{x\in p} (x)^-$ for the negative
topology. A simple argument using the maximality of ultrafilters
shows that these sets are just the complement of $\{p\}$. But since
they are open in the respective topologies, $\{p\}$ is closed in
either topology. \qed\end{pf}

Later on (in Proposition \ref{prop:Hausdorff2} and Proposition
\ref{prop:Hausdorff4}) we shall show that the separation property
${\rm T}_1$ is in general as good as it gets: There are examples of
semilattices $\bbL$ with non-Hausdorff spaces $\YUm{\bbL}$.

For the positive topologies, the property of principality has another
application: In a more special case, we can deduce for principal
spaces that the set of principal ultrafilters is dense in the
positive topology. For this, we shall call a semilattice $\bbL$ {\bf
principally generated} if for each $x\in\bbL$, where $x\neq {\mathbf
0}$, there is a $y\seg x$ such that $[\{y\}]$ is a principal
ultrafilter on $\bbL$. Note that if $\bbL$ is principally generated,
then $\YU{\bbL}$ contains principal ultrafilters.

\begin{obs}\label{obs:dense}
If $\bbL$ is a principally generated semilattice, then the set of
principal ultrafilters is dense in $\YUp{\bbL}$.
\end{obs}

\begin{pf} Let $(x)^+$ be an arbitrary, non-empty basic open set.
By the assumption, there is $y\seg x$ such that $p:=[\{y\}]$ is an
ultrafilter. Thus, $p \in (x)^+$ and hence the set of principal
ultrafilters intersects any open set.\qed\end{pf}

\section{The Ultrafilter Spaces on the Set of Partitions}\label{sec:UltPart}

In order to prove the following results we introduce first some
notation.

\medskip

In the following, for an arbitrary set $x$, let $|x|$ denote the
cardinality of $x$. We always identify a natural number $n\in\omega$
with the set $n=\{m\in\omega: m<n\}$. For $x\subs\omega$ let $\mmin
(x):=\bigcap x$. If $X$ is a partition, then $\mMin (X):=\{\mmin
(x):x\in X\}$; and for $n\in\omega$ and $X\in\parto$, $X(n)$ denotes
the unique block $x\in X$ such that $|\mmin (x)\cap\mMin (X)|=n+1$.
($X(n)$ is just the $n$th block of $X$ in the order of increasing
minimal elements.) Finally, a partition is called {\bf trivial} if it
contains only one block.

Concerning $\YU{\parce}$, we like to mention the following

\begin{fct}\label{fct:fin2}
If $p$ is an ultrafilter on $\parti$ and $p$ contains a finite
partition, then there is a $2$-block partition $X$ such that $p = [\{
X\}]$, and hence, $p$ is principal.
\end{fct}

\begin{pf} Let $m := \min\{ n : \exists Y \in p(|Y|=n)\}$. This minimum exists by
assumption. Let $X\in p$ be such that $|X|=m$.

First we show that for all $Y\in p$ we have $X\ceq Y$. Suppose this
is not the case for some $Y\in p$, then we have $X\neq X\kap Y\in p$
(since $p$ is a filter), which implies $|X\kap Y|<|X|=m$ and
contradicts the definition of $m$. On the other hand, there is a
$2$-block partition $Z$ with $Z\ceq X$, and because $Z\ceq X$ we get
$Z\ceq Y$ for any $Y\in p$. Therefore, since $p$ is an ultrafilter,
we get $Z=X$, which implies $[\{ X\}] = p$ and $m=2$. \qed\end{pf}

This leads to the following observations:

\begin{fct}\label{fct:singletons}
The space $\YUp{\parce}$ is a principal topological space, whereas
the space $\YUp{\parfe}$ is non-principal.
\end{fct}

\begin{pf}
That $\YUp{\parce}$ is principal follows directly from Fact
\ref{fct:fin2} and Proposition \ref{prop:principal}. For the second
assertion we note that for every partition $Y\in\parti$ we find
$Z_1,Z_2\in\parti$ such that $Y\ceq Z_1$, $Y\ceq Z_2$ and $Z_1\kup
Z_2={\mathbf 0}$, and therefore, we find $p_1,p_2\in\YU{\parfe}$ with
$Z_1\in p_1$ and $Z_2\in p_2$, which implies that $p_1$ and $p_2$
both belong to $(Y)^+$. So, for each $Y\in\parti$, the set $(Y)^+$ is
not a singleton. (In fact, by this argument, $\YUp{\parfe}$ doesn't
have any finite open sets.)\qed\end{pf}

\subsection{The space ${\YUp{\parce}}$}

As in the space $\beta\omega$, the principal ultrafilters in
$\YU{\parce}$ form a dense set
in ${\YUp{\parce}}$ by Observation \ref{obs:dense}, but since there
are continuum many 2-block partitions (one for each subset of
$\omega$) in $\YU{\parce}$, they cannot witness that the space
$\YUp{\parce}$ is separable. Moreover, we get the following

\begin{obs}\label{obs:separabel}
The space $\YUp{\parce}$ is not separable.
\end{obs}

\begin{pf} Spinas proved in \cite{Spinas} that there is an uncountable
set $\{X_{\iota}:\iota\in I\}\subs\parto$ of infinite partitions such
that $X_{\iota}\kap X_{\iota'}=\{\omega\}$ whenever
$\iota\neq\iota'$. Thus, $(X_{\iota})^+ \cap (X_{\iota'})^+
=\emptyset$ (for $\iota\neq\iota'$), which implies that there is no
countably dense set in the space $\YUp{\parce}$. \qed\end{pf}

\begin{propn}\label{prop:Hausdorff1} The space $\YUp{\parce}$ is a
Hausdorff space.
\end{propn}

\begin{pf} Let $p$ and $q$ be two distinct ultrafilters $\YU{\parce}$.
Because $p\neq q$ and both are maximal filters, we find partitions
$X\in p$ and $Y\in q$ such that $X\kap Y={\mathbf 0}$. So we get
$p\in (X)^+$, $q\in (Y)^+$ and $(X)^+\cap (Y)^+ =\emptyset$. \qed\end{pf}

Before we prove the next proposition, we state the following useful

\begin{lm}\label{lm:ort}
If $X_0,\ldots,X_n\in\parti$ is a finite set of non-trivial
partitions, then there is a non-trivial partition $Y\in\parti$ such
that $Y\ort_\ceq X_i$ for all $i\le n$.
\end{lm}

\begin{pf} Let $Z_0:=\mMin(X_0)$. If $Z_{i}$ is such that $Z_i\cap
X_{i+1}(k)\neq\emptyset$ for every $k\le |X_{i+1}|$, then
$Z_{i+1}=Z_i$. Otherwise, we define $Z_{i+1}\supseteq Z_i$ as
follows: If $Z_i\cap X_{i+1}(k)\neq\emptyset$, then $Z_{i+1}\cap
X_{i+1}(k)=Z_i\cap X_{i+1}(k)$; and if $Z_i\cap X_{i+1}(k)=
\emptyset$, then $Z_{i+1}\cap X_{i+1}(k)=\mmin(X_{i+1}(k))$. It is
easy to see that $\omega\setminus Z_{i}$ is infinite for every $i\le
n$. Finally, let $Y=\{Y(0),Y(1)\}\in\parti$ be such that $Z_n\subs
Y(0)$ and by construction we get $Y\ort_\ceq X_i$ for all $i\le n$. \qed\end{pf}

\begin{propn}\label{prop:quasicompact1}
The space $\YUp{\parce}$ is not quasicompact.
\end{propn}

\begin{pf} Let $\cA =\{(X)^+ :X\in\parto\}$, then it is easy to see that
$\bigcup\cA =\YU{\parce}$. We will show that A is a cover
with no finite subcovers.
Assume to the contrary that there are
finitely many infinite partitions $X_0,\ldots,X_n\in\parto$ such that
$(X_0)^+\cup\ldots\cup (X_n)^+ =\YU{\parce}$. By Lemma~\ref{lm:ort}
we find a $Y\in\parti$ such that $Y\ort_\ceq X_i$ (for all $i\le n$). Let
$p\in\YU{\parce}$ be such that $Y\in p$, then $X_i\notin p$ (for all
$i\le n$), which contradicts the assumption. \qed\end{pf}

\subsection{The space ${\YUm{\parce}}$}

\begin{propn}\label{prop:Hausdorff2}
The space ${\YUm{\parce}}$ is not a Hausdorff space.
\end{propn}

\begin{pf} Let $p$ and $q$ be two distinct ultrafilters in $\YU{\parce}$.
Take any non-trivial partitions
$X_0,\ldots,X_k,\,Y_0,\ldots,Y_\ell\in\parti$ such that $p\in
(X_0)^-\cap\ldots\cap (X_k)^-$ and $q\in (Y_0)^-\cap\ldots\cap
(Y_\ell)^-$. Now, by Lemma~\ref{lm:ort}, there is a non-trivial
partition $Z$ such that $Z\ort_\ceq X_i$ (for $i\le k$) and $Z\ort_\ceq Y_j$
(for $j\le \ell$), which implies $Z\in \bigcap_{i\le k}(X_i)^- \cap
\bigcap_{j\le \ell}(Y_j)^-$. Hence, $\bigcap_{i\le k}(X_i)^- \cap
\bigcap_{j\le \ell}(Y_j)^-$ is not empty. \qed\end{pf}

\begin{propn}\label{prop:ctblycompact2}
The space ${\YUm{\parce}}$ is countably compact.
\end{propn}

\begin{pf} Let $\cA=\{\bigcap A_i:i\in\omega\}$ be such that $\bigcup\cA=
\bigcup_{i\in\omega} (\bigcap A_i)=\YU{\parce}$, where each $A_i$ is
a finite set of open sets of the form $(X)^-$ for some $X\in\parti$.
Assume $\bigcup_{i\in I}(\bigcap A_i)\neq\YU{\parce}$ for every
finite set $I\subs\omega$. If $A_i=\{(X_0^i)^-,\ldots ,(X_n^i)^- \}$
and $A_j=\{(X_0^j)^-,\ldots, (X_m^j)^- \}$ and $\bigcap A_i\cup
\bigcap A_j\neq\YU{\parce}$, then we find a $p\in\YU{\parce}$ such
that $p\in \YU{\parce}\setminus \bigcap A_i\cup \bigcap A_j$. Hence,
there are $k\le n$ and $\ell\le m$ such that $X_k^i$ and $X_\ell^j$
are both in $p$, which implies $X_k^i\kap X_\ell^j\neq {\mathbf 0}$.
We define a tree $T$ as follows: For $n\in\omega$ the sequence $\la
s_0,\ldots,s_n\ra$ belongs to $T$ if and only if for every $i\le n$
there is an $(X_k^i)^-\in A_i$ such that $s_i=X_k^i$ and
$(s_0\kap\ldots\kap s_n)\neq{\mathbf 0}$. The tree $T$, ordered by
inclusion, is by construction (and by our assumption) a tree of
height $\omega$ and each level of $T$ is finite. Therefore, by
K\"onig's Lemma, the tree $T$ contains an infinite branch. Let $\la
X^i\, : \, i\in\omega\ra$ be an infinite branch of $T$, where $X^i\in
A_i$. By construction of $T$, for every finite
$I=\{\iota_0,\ldots,\iota_n\}\subs\omega$ we have
$X^{\iota_0}\kap\ldots\kap X^{\iota_n}\neq {\mathbf 0}$. Thus the
partitions constituting the branch have the finite intersection
property and therefore we find a $p\in\YU{\parce}$ such that $X^i\in
p$ for every $i\in\omega$. Now, $p\notin\bigcup_{i\in\omega}(X^i)^-$
which implies that $p\notin\bigcup\cA$, but this contradicts
$\bigcup\cA=\YU{\parce}$. \qed\end{pf}

\subsection{The space ${\YUp{\parfe}}$}

\begin{propn}\label{prop:Hausdorff3}
The space $\YUp{\parfe}$ is a Hausdorff space.
\end{propn}

\begin{pf} Let $p$ and $q$ be two distinct ultrafilters $\YU{\parfe}$.
Because $p\neq q$ and both are maximal filters, we find partitions
$X\in p$ and $Y\in q$ such that $X\kup Y={\mathbf 0}$. Hence we get
$p\in (X)^+$, $q\in (Y)^+$ and $(X)^+\cap (Y)^+ =\emptyset$. \qed\end{pf}

Before we prove the next proposition, we state the following useful

\begin{lm}\label{lm:ortfe}
If $X_0,\ldots,X_n\in\parlo$ is a finite set of non-trivial, finite
partitions, then there is a finite partition $Y\in\parlo$ such that
$Y\ort_\feq X_i$ for all $i\le n$.
\end{lm}

\begin{pf} Define an equivalence relation on $\omega$ as follows:
$$s\approx t\; : \iff \forall i,k \big( s\in X_i(k)\leftrightarrow
t\in X_i(k)\big)$$ Because every partition $X_i$ is finite and we
only have finitely many partitions $X_i$, at least one of the
equivalence classes must be infinite, say $I$. Since each block of
each partition $X_i$ is infinite and the partitions have been assumed
to be non-trivial, we also must have $\omega\setminus I$ is infinite.
Let $I_{-1}:=I$ and define $I_{i+1}:=I_i\dot\cup\{s_{i+1}\}$ in such
a way that for any $t\in I$ we have $s_{i+1}\in X_{i+1}(k)\rightarrow
t\notin X_{i+1}(k)$. Let $Y:=\{I_n,\omega\setminus I_n\}$, then
$Y\in\parti$ and for every $i\le n$, $Y\kup X_i$ contains a finite
block and therefore, $Y\ort_\feq X_i$ (for all $i\le n$).
\qed\end{pf}

\begin{propn}\label{prop:quasicompact3}
The space $\YUp{\parfe}$ is not quasicompact.
\end{propn}

\begin{pf} Let $\cA =\{(X)^+ :X\in\parlo\}$, then it is easy to see that
$\bigcup\cA =\YU{\parfe}$. Assume to the contrary that there are
finitely many finite partitions $X_0,\ldots,X_n\in\parlo$ such that
$(X_0)^+\cup\ldots\cup (X_n)^+ =\YU{\parfe}$. By Lemma~\ref{lm:ortfe}
we find a $Y\in\parlo$ such that $Y\ort_\feq X_i$ (for all $i\le n$). Let
$p\in\YU{\parfe}$ be such that $Y\in p$, then $X_i\notin p$ (for all
$i\le n$), which contradicts the assumption. \qed\end{pf}

\subsection{The space ${\YUm{\parfe}}$}

\begin{propn}\label{prop:Hausdorff4}
The space ${\YUm{\parfe}}$ is not a Hausdorff space.
\end{propn}

\begin{pf} We first show that if $p\in (X)^-$ for some $X\in\parto$, then
there is an $X'\in\parlo$ such that $X'\ceq X$ (and therefore
$(X')^-\subs (X)^-$) and $p\in (X')^-$. Since $p\in (X)^-$, there is
a $Y\in p$ such that $Y\kup X={\mathbf 0}$, which is equivalent to
the statement (because we only allowed infinite blocks):
There are $y\in Y$ and $x\in X$ such that $x\cap y$ is
a non-empty, finite set. Now, for $X':=\{x,\omega\setminus x\}$ we
obviously have $X'\ceq X$ and $p\in (X')^-$.

Let $p$ and $q$ be two distinct ultrafilters in $\YU{\parfe}$. Take
any partitions $X_0,\ldots,X_k,\,Y_0,\ldots,Y_l\in\parti$ such that
$p\in (X_0)^-\cap\ldots\cap (X_k)^-$ and $q\in (Y_0)^-\cap\ldots\cap
(Y_l)^-$. By the fact mentioned above we may assume that the $X_i$'s
as well as the $Y_i$'s are finite partitions. Now, by
Lemma~\ref{lm:ortfe}, there is a finite partition $Z$ such that
$Z\ort_\feq X_i$ (for $i\le k$) and $Z\ort_\feq Y_j$ (for $j\le l$), which
implies $Z\in \bigcap_{i\le k}(X_i)^- \cap \bigcap_{j\le l}(Y_j)^-$.
Hence, $\bigcap_{i\le k}(X_i)^- \cap \bigcap_{j\le l}(Y_j)^-$ is not
empty. \qed\end{pf}

\begin{propn}\label{prop:ctblycompact4}
The space ${\YUm{\parfe}}$ is countably compact.
\end{propn}

\begin{pf} Replacing ``$\kap$'' by ``$\kup$'' and ``$\ceq$'' by
``$\feq$'', one can simply copy the proof of
Proposition~\ref{prop:ctblycompact2}. \qed\end{pf}

\subsection{Conclusion} Now we are ready to state the main result of
this paper.

\begin{thrm}\label{thm:main}
None of the spaces $\YUp{\parce}$, $\YUm{\parce}$, $\YUp{\parfe}$ and
$\YUm{\parfe}$ is homeomorphic to $\beta\omega$ or
$\beta\omega\setminus\omega$. Moreover, no two of the spaces
$\beta\omega$, $\beta\omega\setminus\omega$, $\YUp{\parce}$,
$\YUm{\parce}$ and $\YUp{\parfe}$ are homeomorphic.
\end{thrm}

\begin{pf} The proof is given in the following table which is just the
compilation of the results from Sections \ref{sec:UltPart} and
\ref{sec:UltSem}. The separation property ${\rm T}_1$ holds for all
spaces and thus does not help to discern any two spaces; it is just
included for completeness.

{\scriptsize{
\begin{center}
\begin{tabular}{|l||c|c||c|c|c|c|}
\hline & $\beta\omega$ & $\beta\omega\setminus\omega$ & $\YUp{\parce}$ &
$\YUm{\parce}$ & $\YUp{\parfe}$ & $\YUm{\parfe}$ \\[0,1cm] \hline

principal & {\sc Yes} & {\sc No} & {\sc Yes} & & {\sc
No}&\\

${\rm T}_1$ & {\sc Yes} &{\sc Yes} &{\sc Yes} &{\sc Yes} &{\sc Yes} &{\sc Yes}\\

Hausdorff & {\sc Yes} & {\sc Yes} & {\sc Yes} & {\sc No} & {\sc Yes}
& {\sc No}\\

ctb.\;compact & {\sc Yes} & {\sc Yes} & {\sc } & {\sc Yes} & {\sc }&
{\sc Yes}\\

quasicompact & {\sc Yes} & {\sc Yes} & {\sc No} & {\sc } & {\sc No} &
{\sc }\\ \hline
\end{tabular}
\end{center}
}}
{\hfill}\qed\end{pf}

Note that in the language of Section \ref{sec:int}, this immediately implies that
the partial order
$\la\Cosub{\omega},\leq\ra$ of cosubobjects of $\omega$
is not a Boolean algebra (not even a complemented semilattice) since
otherwise we would have $\parce\cong\parfe$ and hence $\YUp{\parce}$ and
$\YUp{\parfe}$ would be homeomorphic.

\section{About the space
$\boldsymbol{\YUp{\partoce}}$}\label{sec:inf}

To investigate the space ${\YUp{\partoce}}$ we first introduce some
notations.

\bigskip

For $X\in\parti$ and $n\in\omega$ let $X\kap\{n\}$ be the partition
we get, if we glue all blocks of $X$ together which contain a member
of $n$. If $X,Y\in\parto$, then we write $X\ceq^* Y$ if there is
an $n\in\omega$ such that $(X\kap \{n\})\ceq Y$. For $X,Y\in\parto$
it is not hard to see that in the space $\YUp{\partoce}$ we have
$(X)^+\subs (Y)^+$ if and only if $X\ceq^* Y$.

\subsection{The height of tree ${\pi}$-bases of
${\YUp{\partoce}}$}

We first give the definition of the dual-shattering cardinal $\fH$.

A family ${\mathcal A}\subseteq\parto$ is called {\bf maximal
orthogonal} (\mo) if ${\mathcal A}$ is a maximal family of pairwise
orthogonal partitions. A family ${\mathcal H}$ of {\mo} families of
partitions {\bf shatters} a partition $X\in\parto$, if there are
$H\in {\mathcal H}$ and two distinct partitions in $H$ which are both
compatible with $X$. A family of {\mo} families of partitions is {\bf
shattering} if it shatters each member of $\parto$. The
dual-shattering cardinal $\fH$ is the least cardinal number $\kappa$,
for which there is a shattering family of cardinality $\kappa$.

The dual-shattering cardinal $\fH$ is a dualization of the well-known
shattering cardinal $\fh$ introduced by Balcar, Pelant and Simon in
\cite{Balcar.etal} where the letter $\fh$ comes from the word
``height''. In \cite{Balcar.etal} it is proved that $$\fh
=\mmin\{\kappa: \text{there is a tree $\pi$-base for
$\beta\omega\setminus\omega$ of $\fh$eight $\kappa$}\}$$ where a family
$\cB$ of non-empty open sets is called a {\bf
$\boldsymbol{\pi}$-base} for a space $S$ provided every non-empty
open set contains a member of $\cB$, and a {\bf tree
$\boldsymbol{\pi}$-base} $T$ is a $\pi$-base which is a tree when
considered as a partially ordered set under reverse inclusion ({\sl
i.e.}, for every $t\in T$ the set $\{s\in T:s\supseteq t\}$ is
well-ordered by $\supseteq$). The height of an element $t\in T$
is the ordinal $\alpha$ such that $\{s\in T:s\supsetneqq t\}$ is of
order type $\alpha$, and the {\bf height} of a tree $T$ is the
smallest ordinal $\alpha$ such that no element of $T$ has height
$\alpha$.

One can show that $\fH\le\fh$ and $\fH\le\fS$, where $\fS$ is the
dual-splitting cardinal (cf.\,\cite{Cichon.etal}).

It is consistent with the axioms of set theory (denoted by
\ZFC) that $\fH=\aleph_2=2^{\aleph_0}$ ({\sl cf.}~\cite{Halbeisen})
and also that $\fH =\aleph_1
<\fh=\aleph_2$ ({\sl cf.}~\cite{Spinas}). Further it is consistent
with \ZFC\,+\,{\sf MA}\,+\,$2^{\aleph_0}=\aleph_2$ that $\fH
=\aleph_1 <\fh =\aleph_2$, where {\sf MA} denotes Martin's Axiom
({\sl cf.}~\cite{Brendle}).

Following Balcar, Pelant and Simon, it is not hard to prove the
following

\begin{propn}\label{prop:H}
Let $\fH$ be the dual-shattering cardinal defined as above, then
$$\fH =\mmin\big{\{}\kappa: \text{there is a tree $\pi$-base for
$\YUp{\partoce}$ of height $\kappa$}\big{\}}\,.$$
\end{propn}

\begin{pf} Having in mind that for every countable decreasing sequence of
basic open sets $(X_0)^+\supseteq (X_1)^+\supseteq\ldots\supseteq
(X_n)^+\supseteq\ldots$ there is a basic open set $(Y)^+$ such that
for all $i\in\omega$ we have $(Y)^+\subs (X_i)^+$ ({\sl
cf.}~\cite[Prop.\,4.2]{Matet}), one can follow the proof of the Base
Matrix Lemma~2.11 of \cite{Balcar.etal}. As a matter of fact we like
to mention that every infinite {\mo} family has the cardinality of
the continuum (cf.\,\cite{Cichon.etal} or \cite{Spinas}). \qed\end{pf}

Because the shattering cardinal and the dual-shattering cardinal
can be different, this gives us an asymmetry between the two spaces
$\beta\omega\setminus\omega$ and $\YUp{\partoce}$.

\subsection{On ${P}$-points in ${\YUp{\partoce}}$}

In this section we give a sketch of the proof that $P$-points exist
under the assumption of the Continuum Hypothesis, and in general
both existence and non-existence of $P$-points are consistent
with the axioms of set theory. To do this, we will use the technique of
forcing ({\sl cf.}~\cite{Jech}).

\bigskip

An ultrafilter $p\in\YUp{\bbL}$ is a {\bf $\boldsymbol{P}$-point} if
the intersection of any family of countably many neighbourhoods of
$p$ is a (not necessarily open) neighbourhood of $p$.

\medskip

First we show that a $P$-point in $\YUp{\partoce}$ induces in a
canonical way a $P$-point in $\beta\omega\setminus\omega$.

\begin{lm}\label{lm:P}
If there is a $P$-point in $\YUp{\partoce}$, then there is a
$P$-point in $\beta\omega\setminus\omega$ as well.
\end{lm}

\begin{pf} Let $p$ be a $P$-point in $\YUp{\partoce}$, then it is not
hard to see that the filter generated by $\{\mMin(X):X\in p\}$ is a
$P$-point in $\beta\omega\setminus\omega$. \qed\end{pf}

\begin{propn}
It is consistent with $\ZFC$ that there are no $P$-points in
$\YUp{\partoce}$.
\end{propn}

\begin{pf} Shelah proved ({\sl cf.}~\cite[Chapter VI, \S 4]{Shelah}) that it
is consistent with $\ZFC$ that there are no $P$-points in
$\beta\omega\setminus\omega$. But in a model
of $\ZFC$ in which there are no $P$-points in
$\beta\omega\setminus\omega$, there are also no $P$-points in
$\YUp{\partoce}$ by Lemma \ref{lm:P}. \qed\end{pf}

Let $\bbU=\la\parto,\le\ra$ be the partial order defined as follows:
$$X\leq Y \ \Leftrightarrow\ X\ceq^* Y\,.$$

The forcing notion $\bbU$ is a natural dualization of
$\cP(\omega)/{\rm fin}$.

\begin{lm}\label{lm:Gp}
If $G_p$ is $\bbU$-generic over ${\bfV}$, then $G_p$ is a $P$-point
in $\YUp{\partoce}$ in the model $\bfV [G_p]$.
\end{lm}

\begin{pf} First notice that the forcing notion $\bbU$ is $\sigma$-closed
({\sl cf.}~\cite[Proposition 4.2]{Matet}) and hence, $\bbU$ does not add
new reals. For every countable set of neighbourhoods $\{N_i\,:\, i\in\omega\}$
of the
filter $G_p$ we find a countable set of
partitions $\{X_i\, : \,i\in\omega\}\subseteq G_p$
such that $(X_i)^+\subs N_i$ and $X_i\ceq^*
X_j$ for $i\ge j$. Now, since every partition $X\in\parto$ can be
encoded by a real number and $\bbU$ does not add new reals, there is
a $\bbU$-condition $Y$ which forces that the sequence
$X_0\,{}^*\!\!\feq X_1\,{}^*\!\!\feq\ldots$ belongs to $\bfV$,
and since $\bbU$ is $\sigma$-closed we find an infinite partition
$Z\ceq Y$ such that $Z\ceq^* X_i$ for every $i\in\omega$. Hence, $Z$
forces that $(Z)^+$ belongs to $\bigcap_{i\in\omega}N_i$ and that $Z$
belongs to $G_p$. \qed\end{pf}

\begin{propn}\label{prop:CH} Assume $\CH$, then there is a
$P$-point in $\YUp{\partoce}$.
\end{propn}

\begin{pf} Assume ${\bfV}\models\CH$. Let $\chi$ be large enough such
that $\cP(\parto)\in H(\chi)$, {\sl i.e.}, the power-set of $\parto$
(in ${\bfV}$) is hereditarily of size $<\chi$. Let $\bfN$ be an
elementary submodel of $\langle H(\chi),\in\rangle$ containing all
the reals of ${\bfV}$ such that $|\bfN|=2^{\aleph_0}$. We consider
the forcing notion $\bbU$ in the model $\bfN$. Since
$|\bfN|=2^{\aleph_0}$, in ${\bfV}$ there is an enumeration
$\{D_\alpha\subs\parto:\alpha<2^{\aleph_0}\}$ of all dense sets of
$\bbU$ which lie in $\bfN$. Since $\bbU$ is $\sigma$-closed and
because ${\bfV}\models\CH$, $\bbU$ is $2^{\aleph_0}$-closed in
${\bfV}$ and therefore we can construct a descending sequence
$\{X_\alpha:\alpha<2^{\aleph_0}\}$ in ${\bfV}$ such that $X_\alpha\in
D_\alpha$ for each $\alpha<2^{\aleph_0}$. Let
$G_p:=\{X\in\parto:X_\alpha\ceq X\ \text{for some $X_\alpha$}\}$,
then $G_p$ is $\bbU$-generic over $\bfN$. By Lemma~\ref{lm:Gp} we
have $\bfN[G_p]\models\text{``there is a $P$-point in
$\YUp{\partoce}$''}$ and because $\bfN$ contains all reals of
${\bfV}$ and every countable descending sequence of basic open sets
$(Y_i)^+$ can be encoded by a real number, the $P$-point $G_p$ in the
model $\bfN[G_p]$ is also a $P$-point in $\YUp{\partoce}$ in the
model ${\bfV}$, which completes the proof. \qed\end{pf}

\section{Open Questions}\label{sec:open}

As we already mentioned, we consider our present work more as a
teaser. This paper leaves many questions open, and one of our goals
is to awaken the interest in delving deeper into this subject matter.
The range of questions reaches from deeper inquiries about the spaces
explored in this paper to completely different spaces derived by the
same general method from semilattices.

We would like to remark that de Groote mentions in \cite{dG} a
non-compactness result similar to Proposition
\ref{prop:quasicompact3} for the Stone space of a Hilbert lattice
$\mathbb{L}(\mathcal{H})$ ({\sl i.e.}, the lattice of closed
subspaces of a given Hilbert space $\mathcal{H}$). This shows that
the general framework (topological properties of ultrafilter spaces
over a semilattice) is of interest for a broader community.

\subsection{Semilattices from Recursion Theory}

Many semilattices with an enormous amount of structure results occur
in Recursion Theory: We could take the semilattice $\la
\mathcal{D},\leq_T\ra$ of Turing degrees ({\sl cf.} \cite{Cooper}),
the semilattice $\la \mathcal{H}, \leq_{\Delta^1_1}\ra$ of hyperdegrees
({\sl cf.} \cite{Hinman}), or any of the multifarious degree
structures derived from reducibility relations on the reals. Note
that none of these structures can be complemented semilattices: We
mentioned that in any complemented semilattice $\la L,\seg\ra$ the
semilattice is homomorphic to $\la L,\succcurlyeq\ra$. But for all
degree structures derived from reducibility relations we know that
the set $\{ d : d\seg e\}$ is countable for each $e$, whereas in most
cases the sets $\{d : d\succcurlyeq e\}$ are uncountable.

Therefore the investigation of the ultrafilter spaces on these
semilattices seems to be promising. In addition to these full degree
structures, Recursion Theory has to offer the semilattice of
recursively enumerable sets ({\sl cf.} \cite{Soare}) and similar
structures. For a general overview about known results from Recursion
Theory, we refer the reader to \cite{Slaman}.

\subsection{Other Partition Semilattices}

But we don't have to leave the realm of partitions to find new
mathematically interesting semilattices. In Section \ref{sec:inf}, we
already started to inquire into the properties of the space
$\YUp{\partoce}$. This space can be considered as the dualization of
$\beta\omega\setminus\omega$ and its properties are similar to those
of the space $\YUp{\parfe}$. Because $\kap$ and $\kup$ are not
inverse to each other, it is unlikely that the spaces
$\YUp{\partoce}$ and $\YUp{\parfe}$ are homeomorphic. Further, a
dual-shattering cardinal $\fH^*$ can also be defined in the space
$\YUp{\parfe}$ ({\sl cf.}~\cite{Cichon.etal}). What is the relation
between $\fH^*$ and $\fH$ and between $\fH^*$ and $\fh$?

\bigskip

We had restricted our attention
to partitions which have only infinite blocks. What happens if we
consider {\em all\/} partitions of $\omega$? Let $\parti^{\st}$
denote the set of all possible partitions of $\omega$. What are the
topological properties of the spaces $\YUp{\parce^{\st}}$,
$\YUm{\parce^{\st}}$, $\YUp{\parfe^{\st}}$ and $\YUm{\parfe^{\st}}$?
Are they all different and is one of them homeomorphic to $\beta\omega$
or $\beta\omega\setminus\omega$? What is the relation (if there is any)
between $\YUp{\parce}$ and $\YUp{\parce^{\st}}$ (and likewise for all other pairs)?

We can boldly step forward along the ordinals and look at partitions
of larger sets than $\omega$. We could compare spaces of ultrafilters
of partitions on $\omega_1$ with our spaces. This comparison is to be
seen in connection with the open question whether
$\beta\omega\setminus\omega$ can be homeomorphic to
$\beta\omega_1\setminus\omega_1$ ({\sl cf.} \cite[\S 5, Problem
1]{vanMill}).

\subsection{Deeper Knowledge about the Spaces presently under Investigation}

After listing more natural candidates for the underlying semilattice,
we now proceed to ask deeper questions about the four spaces
scrutinized in this paper.

In this paper our main focus was to show that the spaces are different
from each other and from $\beta\omega$ and
$\beta\omega\setminus\omega$. One small question is left open by
Theorem \ref{thm:main}, though: Is there any topological property
which we could use to distinguish between $\YUm{\parce}$ and
$\YUm{\parfe}$? Regardless of the answer to this more technical
question, the non-homeomorphicity results are supposed to be a
starting point from which one could now move on to derive more
properties of the spaces.

At first, one could investigate properties related to the ones that
we checked in this paper ({\sl e.g.}, separability or the Lindel\"of
property).
For the noncompact positive spaces $\YUp{\parce}$ and $\YUp{\parfe}$ it is not even
obvious that they are compactifiable. To this end, we'd have to show that they are
${\rm T}_{3{\rm a}}$ ({\sl cf.} \cite[Theorem 3.5.1]{Engelking}).
Then, and probably more interesting, one could go back to
the original motivation: the spaces of partition ultrafilters were
constructed as a dualization of the Stone-{\v C}ech compactification
of the natural numbers. From this point on, one could try to find
similarities and differences between our spaces and the Stone-{\v
C}ech compactification.

To name but a few:
\begin{enumerate}
\item One of the most important characterizations of the Stone-{\v C}ech
compacti\-fication uses the {\bf extension of mappings} concept:

$\beta\omega$ is the unique space $X$ such that every continuous map
from $\omega$ to a compact space can be uniquely extended to a
continuous map defined on $X$.

\begin{center}
\begin{minipage}{4cm}
\xymatrix{ \omega \ar[ddrr]^{f}\ar[rr]^{\iota}&
&X=\beta\omega\ar[dd]^{\hat f} \\ &&\\ &&K }
\end{minipage}
\end{center}

Can we show or refute anything analogous for any of our spaces of
partition ultrafilters?

This question seems to be connected to the categorial dualization
process by which we moved from subsets to partitions in Section
\ref{sec:int}. It could lead to more inquiries by looking at other
categoral characterizations and properties of $\beta\omega$ ({\sl
cf.} \cite[\S 10]{Walker}).

\item Building on the previous point, if
we have some extension principle for continuous maps, then we could
introduce the notion of {\bf ultrafilter types} as usual for
$\beta\omega$ ({\sl cf.} \cite[3.41]{Walker}). What can be said about
types of ultrafilters in this sense over partition semilattices? Even
if the extension of mappings does not work, a classification of the
points of the ultrafilter spaces according to some measure of
complexity of the underlying partitions should be possible.

\item Deeply connected with ultrafilter types is the {\bf Rudin-Keisler
order} of ultrafilters ({\sl cf.} \cite[\S 9]{CoNe}. In this field we
would also expect that independently from the success or failure to
get a extension of mappings principle, we should be able to stratify
the ultrafilters according to complexity. Note that the Rudin-Keisler
order provides us with a possibility of characterizing special
ultrafilters on $\omega$. Could there be a topological description of
the Ramsey${}^{{\scriptscriptstyle{\bigstar}}}$ ultrafilters of
\cite{HaLoe} via a dualized Rudin-Keisler ordering?

\item Especially interesting seems the question of autohomeomorphisms
of ultrafilter spaces. For $\beta\omega$ we know that under {\sf CH}
there are many autoho\-meo\-morphisms (\cite[Lemma 1.6.1]{vanMill}), but
Shelah has shown ({\sl cf.} \cite[Section 2.6]{vanMill}) that
consistently, every autohomeomorphism is induced by a permutation of
$\omega$. Results like this become a different feel when we work in
ultrafilter spaces which are not separable like $\YUp{\parce}$.

\item Basically, one can take any part of the abundant literature on
$\beta\omega$ and try to understand the dualized version.
\end{enumerate}

We feel that we have illustrated how little is known about the
fascinating field of ultrafilters on semilattices, and we would like
to see more results along these lines.

\end{document}